\documentclass[english]{article}
\usepackage{amssymb}
\usepackage{amsmath}
\usepackage{babel}

\def\thebibliography#1{\section*{References}\list
  {[\arabic{enumi}]}{\settowidth\labelwidth{[#1]}\leftmargin\labelwidth
    \advance\leftmargin\labelsep
    \usecounter{enumi}}
    \def\newblock{\hskip .11em plus .33em minus -.07em}
    \sloppy
    \sfcode`\.=1000\relax}

\newcommand{\refbook}[3]{{\sc #1}{\em\ #2}{\ #3}}
\newcommand{\refer}[5]{{\sc #1}{\ #2}{\em\ #3}{\bf\ #4}{\ #5}}

\newtheorem{lem}{Lemma}[section]
\newtheorem{cor}[lem]{Corollary}
\newtheorem{teo}[lem]{Theorem}

\newtheorem{defi}[lem]{Definition}
\newtheorem{prop}[lem]{Proposition}

\newcommand{\qed}{\thinspace\null\nobreak\hfill\hbox{\vbox{\kern-.2pt\hrule
 height.2pt depth.2pt\kern-.2pt\kern-.2pt \hbox to2.5mm{\kern-.2pt\vrule
 width.4pt \kern-.2pt\raise2.5mm\vbox to.2pt{}\lower0pt\vtop
 to.2pt{}\hfil\kern-.2pt \vrule
 width.4pt \kern-.2pt}\kern-.2pt\kern-.2pt\hrule height.2pt depth.2pt
 \kern-.2pt}}\par\medbreak}

\newcommand{\R}{\mathbb{R}}

\newcommand{\ds}{\displaystyle}

\date{}
\begin{document}
\title{Optimal kernel estimates for a Schr\"odinger type operator} 
\author{Anna Canale, Cristian Tacelli}
\maketitle
\begin{abstract}
In the paper the principal result obtained is the estimate
for the heat kernel associated to the Schr\"odinger type operator
$(1+|x|^\alpha)\Delta-|x|^\beta$ 
\[
k(t,x,y)\leq Ct^{-\frac{\theta}{2}}\frac
{\varphi(x)\varphi(y)}{1+|x|^\alpha},
\]
where
$\varphi=(1+|x|^\alpha)^{\frac{2-\theta}{4}+\frac{1}{\alpha}\frac{\theta-N}{2}}$,
$\theta\geq N$ and $0<t<1$, provided that
$N>2$, $\alpha> 2$ and $\beta>\alpha-2$. 
This estimate improves a similar estimate in \cite {can-rhan-tac2} with respect 
to the dependence on spatial component.

\par


\end{abstract}

\bigskip

\section{Introduction}\label{sc:quadratic-form}
In this paper we consider the elliptic operator defined by

\begin{equation}\label{eq:operatore-A}
 Au(x)=a(x)\Delta u(x)-V(x)u(x),\qquad\quad x\in\R^N\,,
\end{equation}
where $a(x)=1+|x|^\alpha$, $\alpha>2$ and $V(x)=|x|^\beta$, $\beta>\alpha-2$. 
Our aim is to give better estimate for the associated heat kernel
then those obtained in \cite{can-rhan-tac2}.


Recently elliptic operators with unbounded 
coefficients have been studied in several paper
(see for example \cite{met-spi2}, \cite{met-spi3}, \cite{met-spi-tac}, \cite{met-spi-tac2},
\cite{lor-rhan}, \cite{for-lor}, \cite{can-rhan-tac1}, \cite{kunze-luc-rha1}, \cite{kunze-luc-rha2},
\cite{can-tac1}, \cite{dur-man-tac1}).
In \cite{lor-rhan}  and  \cite{can-rhan-tac1}
it is proved that 
$A$ endowed  with domain
\begin{equation}\label{eq:domain}
D_p(A)=\{u\in W^{2,p}(\R^N)\;|\; (1+|x|^\alpha)|D^2u|,(1+|x|^\alpha)^{1/2}\nabla u,|x|^\beta u\in L^p(\R^N) \} 
\end{equation} 
generates a strongly continuous and analytic semigroup $T(\cdot)$ 
in $L^p(\R^N)$ for $1<p<\infty$,
for $\alpha>2$ and $\beta >\alpha -2$.
This semigroup is also consistent, irreducible and ultracontractive. 
As regards the case $\beta=0$ we refer to \cite{for-lor} and \cite{met-spi2}.


Due to the regularity of the coefficients of the operator $A$,
the semigroup $T(t)$ 
can be represented in the following integral form through a heat kernel $k(t,x,y)$ 
\[
 T(t)f(x)=\int_{\R^N} k(t,x,y)f(y)dy\;,\qquad t>0,\;x\in \R^N\;,
\]
for any $f\in L^p$ (see \cite{ber-lor},\cite{met-wack}).

In \cite{can-rhan-tac2} was obtained the heat kernel estimate
provided that $N>2,\,\alpha\geq 2$ and $\beta>\alpha-2$

\begin{equation}\label{eq:stima-con-autovvalori}
k(t,x,y)\leq c_1 e^{\lambda_0 t}e^{c_2 t^{-b}}
\frac{\psi(x)\psi(y)}{1+|y|^\alpha}, \quad t>0,\,x,\,y\in \R^N, 
\end{equation} 
where $c_1,c_2$ are positive constant, $b=\frac{\beta-\alpha+2}{\beta+\alpha-2}$
and $\psi(x)$ is the eigenfunction associated to the first eigenvalue, which is equivalent to the function
$$|x|^{-\frac{N-1}{2}-\frac{\beta-\alpha}{4}} 
e^{-\int_1^{|x|}\frac{s^{\beta/2}}{\sqrt{1+s^\alpha}}\,ds}.$$
A better estimate with respect to the time variable 
was also obtained for small values of $t$
\begin{align}\label{eq:stima-old-2}
k(t,x,y)\leq Ct^{-\frac{N}{2}}\left(1+|x|^\alpha \right)^{\frac{2-N}{4}}
\left(1+|y|^\alpha \right)^{\frac{2-N}{4}-1}\,,\quad 0<t\leq 1.
\end{align}

Comparing \eqref{eq:stima-con-autovvalori} and \eqref{eq:stima-old-2} ones can see that
improving the dependence on $t$ 
involves worsening in the dependence on the space component. 
Conversely, if the space component is improved the other worsens. 

In this paper our aim is to explain how this happens.
In order to state the relationship between the dependence on 
the time and the space components we will state an estimate which depend 
on a parameter $\theta$.
In particular we will prove the following estimate for small values of $t$
\begin{equation}\label{eq:stima-kernel}
k(t,x,y)\leq Ct^{-\frac{\theta}{2}}\frac
{\varphi(x)\varphi(y)}{1+|y|^\alpha}, 
\end{equation} 
where
$\varphi=(1+|x|^\alpha)^{\frac{2-\theta}{4}+\frac{1}{\alpha}\frac{\theta-N}{2}}$,
$\theta \geq N$.
We observe that \eqref{eq:stima-old-2} is a particular case of \eqref{eq:stima-kernel} 
obtained for $\theta=N$.
\bigskip

\section{Weighted spaces and Weighted Nash Inequalities}

Let $T(t)$ be the semigroup generated by the operator $(A,D_{p}(A))$, where 
$A$ and $D_{p}(A)$ are defined by \eqref{eq:operatore-A} and \eqref{eq:domain}.
First we show that
$T(t)$ can be seen as a suitable  semigroup ${\cal T}(t)$ on a weighted space.
So, we can deduce heat kernel estimates of $T(t)$ by heat kernel estimates of ${\cal T}(t)$.
Then, in order to obtain kernel estimates of ${\cal T}(t)$ 
we prove a weighted ultracontractivity of the semigroup obtained by a
weighted Nash Inequality.

Let us introduce the measure $d\mu(x)=\left 
(1+|x|^\alpha\right)^{-1}dx$ and the Hilbert space 
$L^2_\mu$  endowed with its canonical inner product.
Let  
$$
H=\{u \in L^2_\mu\cap W^{1,2}_{loc} :\; V^{1/2}u\in L^2_\mu\;,\; \nabla u \in L^2\}
$$
be the Sobolev space endowed with the inner product
$$
(u,v)_H=\int_{\R^N} (1+V)u\bar{v}\, d\mu +\int_{\R^N}\nabla u \cdot \nabla
 \bar{v}\, dx. 
$$
We consider the close and accretive  symmetric form so defined
\begin{equation} \label{formaquad}
a(u,v)= \int_{\R^N}\nabla u \cdot \nabla
 \bar{v}\, dx+\int_{\R^N} u \bar v d\mu
\end{equation}
for $u,v$ belonging to the closure $\cal V$ of $C_c^\infty$ in $H$,
with respect to the norm of $H$. 
Then we can associate to $a$ the self-adjoint operator 
$$
{\cal A}u=f
$$ 
with domain
$$
 D({\cal A})=\{u\in {\cal V} : {\rm there\ exists\ } f \in L^2_\mu
\quad s.t. \quad
 a(u,v)=-\int_{\R^N}f\bar{v}\, d\mu {\rm \ for\ any\ } v \in {\cal V} \}.
$$
By classical results the operator ${\cal A}$ generates an analytic semigroup 
of contractions ${\cal T}(t)$ in $L^2_\mu$ 
which is a positive, symmetric and $L^\infty$-contractive Markov semigroup. 

The Lemma below (see \cite{can-rhan-tac2}) shows that the semigroup 
${\cal T}(t)$ coincides in $L^p \cap L^2_\mu$
with the 
semigroup $T(t)$ generated by  $(A,D_{p}(A))$ in $L^p(\R^N)$.
\begin{lem}\label{coerenza1} 
We get
$$
D({\cal A})\subset \left\{u \in {\cal V} \cap W^{2,2}_{loc}:
(1+|x|^\alpha) \Delta u -V(x)u \in L^2_\mu \right\}
$$
and ${\cal A} u=(1+|x|^\alpha)\Delta u-V(x)u $ for $u \in D({\cal A})$. If
$\lambda
>0$ and $f \in L^p \cap L^2_\mu$, then
$$
(\lambda-{\cal A})^{-1}f=(\lambda-A)^{-1}f.
$$
\end{lem}

Denoting  by $k(t,x,y)$ and $k_\mu(t,x,y)$
respectively the heat kernel associated to $T(t)$ and ${\cal T}(t)$,
by the previous lemma we can deduce that
\[
 k_\mu(t,x,y)=(1+|y|^\alpha)k(t,x,y).
\]

In the following we describe how estimates of the kernel
of a symmetric Markov semigroup
can be obtained
by using the equivalence between
a weighted Nash inequality and  a ``weighted''
ultracontractivity of the semigroup.
The equivalence was stated  in \cite[Theorem 3.3]{wang} 
and  was reformulated in \cite[Theorem 2.5]{bakry}.
The equivalence is obtained by means of 
a suitable Lyapunov functions for the generator of the semigroup.

Let ${\cal T}(t)$ be a symmetric Markov semigroup  generated by a self-adjoint operator $\cal A$ associated to 
an accretive, closed, symmetric form $a$ defined on a domain $\cal V$ in $L^2_\mu$.
Let $k_\mu$ its associated heat kernel. We define Lyapunov function in the following way (see also \cite{met-spi3},\cite{met-spi-tac})

\begin{defi} \label{lyap}
 A Lyapunov function is a positive function $\varphi \in L^2_\mu$   such that
$${\cal T}(t)\varphi(x)=\int_{\R^N}k_\mu(t,x,y) \varphi(y)d \mu (y)\leq e^{\kappa t}\varphi(x)$$
for any $x\in\R^N$, $t>0$, and for some real constant $\kappa$, called Lyapunov constant.
\end{defi}





Now, we define the weighted Nash inequality.

\begin{defi}  \label{weightedNash}
Let $\varphi$ be a positive function on $\R^N$  and $\psi$ be a positive
function defined on $(0,\infty)$ with $\ds\frac{\psi(x)}{x}$ non decreasing.
The form $a$ on $L^2_\mu$ satisfies a weighted
Nash inequality with weight $\varphi$ and rate function $\psi$ if 
$$
\psi\left(\frac{\|u\|^2_{L^2_\mu}}{\|u\varphi\|^2_{L^1_\mu}}\right)
\leq\frac{a(u,u)}{\|u\varphi\|^2_{L^1_\mu}}
$$
for any functions $u \in \cal V$ such that 
$\|u\|^2_{L^2_\mu}>0$ and $\|u\varphi\|^2_{L^1_\mu}<\infty$.
\end{defi}

\begin{teo}\cite[Theorem 2.5]{bakry} \label{bakry}
Let ${\cal T}(t)$ be a Markov semigroup with generator 
${\cal A}$ symmetric in $L^2_{\mu}(\R^N)$. 
Let us assume that there exists a Lyapunov function $\varphi$
with Lyapunov constant $\kappa\geq 0$ and that the associated  form $a$
satisfies a weighted Nash inequality with weight
$\varphi$  and rate function $\psi$
 such that 
\[
 \int^\infty \frac{1}{\psi(x)}dx<\infty.
\]
Then
$$
\|{\cal T}(t)f\|_{L^2_\mu}\leq K(2t)e^{\kappa t}\|f\varphi\|_{L^1_\mu}
$$ 
for any functions $f\in L^2_\mu$ such that $f\varphi \in L^1_\mu$.
The function $K$ is defined by
\begin{displaymath}
K(t)=\sqrt{U^{-1}(t)},
\end{displaymath}
 where
$$U(t)=\int_t^\infty\frac{1}{\psi(u)}du.$$
\end{teo}

\noindent Finally from \cite[Corollary 2.8]{bakry} we get the  estimate of $k_\mu$
\begin{cor}\label{cor:pmu}
If the Markov semigroup ${\cal T}(t)$ satisfies the assumptions  
of Theorem \ref{bakry} then the kernel $k_\mu$ satisfies
$$
k_\mu(2t,x,y)\leq K(2t)^2e^{ 2\kappa t}\varphi(x)\varphi(y)
$$ 
for any $t>0$, $(x,y)\in\R^N\times\R^N$.
\end{cor}

\bigskip

\section {Heat kernel estimates}

In this section we will prove upper bound estimates for the kernel $k$.
First we prove that the function $\varphi(x)=
\left(1+|x|^\alpha \right)^\frac{\gamma}{\alpha}$ is a Lyapunov function
if $\gamma<-\frac{N}{2}+\frac{\alpha}{2}$.

\begin{lem}  \label{lyapunov-general}
Let $\gamma<-\frac{N}{2}+\frac{\alpha}{2}$ be a real constant. 
Then the function $\varphi (x)=(1+|x|^\alpha)^{\frac{\gamma}{\alpha} }
\in L^2_\mu(\R^N)$ and satisfies 
the inequality $A\varphi(x)\leq \kappa \varphi $ for some $\kappa>0$.
\end{lem}

{\sc Proof.}
Let us consider $\varphi (x)=(1+|x|^\alpha)^{\frac{\gamma}{\alpha}}$. 
It is easy to see that $\varphi \in L^2_\mu$ if $2\gamma-\alpha<-N$.
Furthermore we get

\begin{align*}
&A\varphi =\gamma \left( \gamma-\alpha\right) 
\left( 1+|x|^\alpha \right)^{\frac{\gamma}{\alpha} -1}|x|^{2\alpha-2}
 +\gamma(\alpha-2+N) \left( 1+|x|^\alpha \right)^{\frac{\gamma}{\alpha}}|x|^{\alpha-2}
 -|x|^\beta \left(1+|x|^\alpha\right)^{\frac{\gamma}{\alpha}}\\
&\quad = \gamma \left( \gamma-\alpha\right) \frac{|x|^{2\alpha-2}}{1+|x|^\alpha}\varphi(x)
 +\gamma(\alpha-2+N) |x|^{\alpha-2}\varphi(x)
 -|x|^\beta \varphi(x)\\
&\quad = \left[
|x|^{\alpha-2}\left( 
\gamma\left( \gamma-\alpha\right) \frac{|x|^{\alpha}}{1+|x|^\alpha}
 +\gamma(\alpha-2+N) \right)
 -|x|^\beta \right]\varphi(x). 
\end{align*}


Then, since $\beta>\alpha-2$, one can see that there exists a positive constant $k$ such that
\[
 A\varphi(x) \leq k\varphi(x) 
\]

\qed

Arguing as in \cite[Section 2]{met-spi2} we have that $\varphi $ is actually a Lyapunov function.

\begin{teo}
The function $\varphi$ is a Lyapunov function with constant $\kappa_0$ 
for any $\kappa_0>\kappa$.
\end{teo}

{\sc Proof.}
Let us observe that $\varphi\in C_0(\R^N)$. Then we can consider 
$u=R(\lambda,A_{\min})\varphi=(\lambda-A_{\min})^{-1}\varphi\in C_0(\R^N)$
(see \cite[Section 2]{can-rhan-tac1}).
Let $\kappa_0>\kappa$, $\lambda\geq \frac{\kappa \kappa_0}{\kappa_0-\kappa}$ and 
$w=\left(1+\frac{\kappa_0}{\lambda}\right)\varphi-\lambda u$.
Since  $Au=\lambda u-\varphi$, we have 

\begin{align*}
&Aw-\lambda w=\frac{\kappa_0+\lambda}{\lambda}A\varphi -\lambda (\lambda u-\varphi)-
\lambda\left(\frac{\kappa_0+\lambda}{\lambda}\varphi-\lambda u  \right)\\
&\quad =\frac{\kappa_0+\lambda}{\lambda}A\varphi -\lambda ^2 u+\lambda\varphi-
(\kappa_0+\lambda)\varphi+\lambda^2 u\\
&\qquad \leq \frac{\kappa_0+\lambda}{\lambda}\kappa \varphi -\kappa_0\varphi=
\frac{1}{\lambda}(\kappa_0\kappa+\lambda\kappa-\lambda \kappa_0)\varphi\leq 0.
\end{align*}
By the maximum principle we have $w>0$ in $\R^N$. Then 
\[
\left( 1+\frac{\kappa_0}{\lambda} \right)\varphi\geq \lambda R(\lambda,A)\varphi.
\]
Iterating the last inequality we get
\[
\left( 1+\frac{\kappa_0}{\lambda} \right)^n\varphi\geq \lambda^n R^n(\lambda,A)\varphi.
\]
So, we obtain
\[
 T(t)\varphi=\lim_{n\to \infty}\left[ \frac{n}{t} R\left(\frac{n}{t},A\right)\right]^n\varphi
 \leq \lim_{n\to \infty}\left( 1+\frac{\kappa_0 t}{n} \right)^n\varphi =e^{\kappa_0 t}\varphi.
\]
\qed

\bigskip

Finally to get kernel estimates we will prove the weighted Nash inequality 
(see Definition \ref{weightedNash}) with Lyapunov function 
$$
\varphi=(1+|x|^\alpha)^{\frac{2-\theta}{4}+\frac{1}{\alpha}\frac{\theta-N}{2}}
$$
and rate functions 
$$
\psi(t)=t^{1+\frac{2}{\theta}}
$$
for $\theta\geq N$. We observe that $\varphi$ satisfies hypothesis Lemma \ref{lyapunov-general} if $\alpha>2$.   

In order to prove the weighted Nash inequality we use a weighted Sobolev inequality which
we recall for reader's convenience (see \cite[Proposition 3.5]{met-spi3}).
\begin{prop} \label{weighted Sobolev} 
Let  $\beta',\ \gamma',\ \nu, p, q$ real values such that 
$$
1< p \le q<\infty
\qquad
\gamma'-1\leq \beta'\leq\gamma',
$$
$$
0\leq\frac{1}{p}-\frac{1}{q}=\frac{1-\gamma'+\beta'}{N}, \quad\quad
N+p(\gamma'-1)\neq 0, \quad \quad p\leq q\leq p^*, \quad\quad p<N.
$$
Then there exists a positive constant $C$ such that for any $u\in
C_c^\infty(\R^N)$
\begin{align*}\label{weighted-Sobolev}
&\left(\int_{\R^N}(1+|x|)^{q\beta'}|u(x)|^qdx\right)^\frac{1}{q}\leq C
\left(\int_{\R^N}(1+|x|)^{\gamma' p}|\nabla u(x)|^pdx\right)^\frac{1}{p}\nonumber 
\\
&\qquad + C\left(\int_{\R^N}(1+|x|)^{\nu}|u(x)|^pdx\right)^\frac{1}{p}. 
\end{align*}
\end{prop}

\begin{teo}  \label{nucleo-tpiccoli}
If $\alpha>2$ and $\beta>\alpha-2$, then the  kernel $k_\mu$ of the semigroup
generated by $A$ satisfies the inequality
\begin{align*}
k_\mu(t,x,y)\leq \frac{C}{t^\frac{\theta}{2}}\varphi(x)\varphi(y)
\end{align*}
 for every $0<t\leq 1$, $x,\ y\in\R^N$.
\end{teo}{\sc Proof.} Let $u\in{\cal V}$ such that
$\|u\varphi\|_{L^1_\mu}<\infty$.

Applying  H\"{o}lder's inequality with $p=\frac{\theta+2}{\theta-2}$ we get
\begin{align*}
&	\int_{\R^N}|u|^2 d\mu
  = \int_{\R^N}|u|^{\frac{2\theta}{\theta+2}+ \frac{4}{\theta+2}} 
\frac{\varphi^\frac{4}{\theta+2}}{\varphi^\frac{4}{\theta+2}}d\mu 
\\
&\qquad	\leq\left(\int_{\R^N}|u|^\frac{2\theta}{\theta-2} 
\varphi^\frac{4}{2-\theta}d\mu\right)^\frac{\theta-2}{\theta+2} 
\left(\int_{\R^N}|u| \varphi d\mu\right)^\frac{4}{\theta+2}
\\
&\qquad	\leq\left(\int_{\R^N}|u|^\frac{2\theta}{\theta-2} 
    (1+|x|^\alpha)^
      {\frac{4}{2-\theta} \left( {\frac{2-\theta}{4}+
\frac{1}{\alpha}\frac{\theta-N}{2}}\right)-1} dx\right)^\frac{\theta-2}{\theta+2}
    \left(\int_{\R^N}|u| \varphi d\mu\right)^\frac{4}{\theta+2}
\\
&\qquad	\leq\left(\int_{\R^N}|u|^\frac{2\theta}{\theta-2} 
    (1+|x|)^
      {2\frac{\theta-N}{2-\theta}} dx\right)^\frac{\theta-2}{\theta+2}
    \left(\int_{\R^N}|u| \varphi d\mu\right)^\frac{4}{\theta+2}.
\end{align*}
Then
\[
 \left( \| u\|^2_{L^2_\mu}  \right)^{1+\frac{2}{\theta}}\leq 
 \left(\int_{\R^N}|u|^\frac{2\theta}{\theta-2} 
    (1+|x|)^
      {2\frac{\theta-N}{2-\theta}} dx\right)^\frac{\theta-2}{\theta}
 \|u\varphi\|_{L^1_\mu}^{\frac{4}{\theta}},
\]
from which
\[
\psi\left(\frac{\|u\|^2_{L^2_\mu}}{\|u\varphi\|^2_{L^1_\mu}}\right)\leq
\frac
  {\left(\int_{\R^N}|u|^\frac{2\theta}{\theta-2} 
    (1+|x|)^
      {2\frac{\theta-N}{2-\theta}} dx\right)^\frac{\theta-2}{\theta}}
  {\|u\varphi\|^2_{L^1_\mu}}
\]
Applying the weighted Sobolev inequality
with $p=2$, $q=\frac{2\theta}{\theta-2}$, 
$\gamma'=0$, $q\beta'=2\frac{\theta-N}{2-\theta}$ and $\nu=-\alpha$
we obtain 
\[
\left(\int_{\R^N}|u|^\frac{2\theta}{\theta-2} 
    (1+|x|)^
      {2\frac{\theta-N}{2-\theta}} dx\right)^\frac{\theta-2}{\theta}
      \leq C\left( \|\nabla u\|^2_2+\int_{\R^N}(1+V)u^2d\mu\right)=C\tilde a(u,u),
\]
where $\tilde a(u,u)=a(u,u)+\int u^2d\mu$ is the quadratic form associated with the operator
$A+I$. Since $\varphi$ is a Lyapunov function with constant $\kappa+1$ for the operator $A+I$,
applying Corollary \ref{cor:pmu}, we get
\[
\tilde k_\mu(t,x,y)\leq \frac{Ce^{(\kappa +1)t}}{t^\frac{\theta}{2}}\varphi(x)\varphi(y), 
\]
where $\tilde k_\mu=e^{t}k_\mu$ is the kernel associated with $A+I$.
This gives the result. 
\qed

Taking into account the relation between $k$ and $k_{\mu}$,
by Theorem \ref{nucleo-tpiccoli} 
we state 
the following result.

\begin{teo}
Le us assume $\alpha >2$, $\beta>\alpha-2$. Then the  kernel $k$ 
of the semigroup generated by $A$ for every $0\leq t \leq 1$ satisfies the bound
\begin{equation}
k(t,x,y)\leq C\frac{e^{\kappa t}}{t^\frac{\theta}{2}}\frac
{\varphi(x)\varphi(y)}{1+|x|^\alpha}.
\end{equation}
\end{teo}

\bigskip

\bibliography{bibfile}{}

\providecommand{\bysame}{\leavevmode\hbox to3em{\hrulefill}\thinspace}
\providecommand{\MR}{\relax\ifhmode\unskip\space\fi MR }
\providecommand{\MRhref}[2]{%
  \href{http://www.ams.org/mathscinet-getitem?mr=#1}{#2}
}
\providecommand{\href}[2]{#2}
\begin{thebibliography}{10}

\bibitem{bakry}
D.~Bakry, F.~Bolley, I.~Gentil, and P.~Maheux, \emph{Weighted nash
  inequalities}, arXiv: 1004.3456.

\bibitem{ber-lor}
M.~Bertoldi and L.~Lorenzi, \emph{Analytical methods for markov semigroups},
  Chapman \& Hall/CRC, 2007.

\bibitem{can-rhan-tac2}
A.~Canale, A.~Rhandi, and C.~Tacelli, \emph{Kernel estimates for
  schr{\"o}dinger operators with unbounded diffusion and potential terms},
  Preprint.

\bibitem{can-rhan-tac1}
\bysame, \emph{Schr{\"o}dinger type operators with unbounded diffusion and
  potential terms}, to appear in Ann. Sc. Norm. Super. Pisa Cl. Sci.

\bibitem{can-tac1}
A.~Canale and C.~Tacelli, \emph{Optimal kernel estimates for a schr{\"o}dinger
  type operator}, Preprint.

\bibitem{dur-man-tac1}
T.~Durante, R.~Manzo, and C.~Tacelli, \emph{Kernel estimates for
  schr{\"o}dinger type operators with unbounded coefficients and singular
  potential terms}, Preprint.

\bibitem{for-lor}
S.~Fornaro and L.~Lorenzi, \emph{Generation results for elliptic operators with
  unbounded diffusion coefficients in ${L}^p$- and ${C}_b$-spaces}, Discr.
  Cont. Dyn. Syst. A \textbf{18} (2007), 747--772.

\bibitem{kunze-luc-rha2}
M.~Kunze, L.~Lorenzi, and A.~Rhandi, \emph{Kernel estimates for nonautonomous
  kolmogorov equations with potential term}, New prospects in direct, inverse
  and control problems for evolution equations \textbf{Springer INdAM Ser., 10}
  (2014), 229--251.

\bibitem{kunze-luc-rha1}
\bysame, \emph{Kernel estimates for nonautonomous kolmogorov equations}, Adv.
  Math. \textbf{287} (2016), 600--639.

\bibitem{lor-rhan}
L.~Lorenzi and A.~Rhandi, \emph{On {S}chr{\"o}dinger type operators with
  unbounded coefficients: generation and heat kernel estimates}, J. Evol. Equ.
  \textbf{15} (2015), 53--88.

\bibitem{met-wack}
G.~Metafune, D.~Pallara, and M.~Wacker, \emph{Feller semigroups on
  $\mathbb{R}^n$}, Semigroup Forum \textbf{65} (2002), 159--205.

\bibitem{met-spi2}
G.~Metafune and C.~Spina, \emph{Elliptic operators with unbounded coefficients
  in ${L}^p$ spaces}, Ann. Sc. Norm. Super. Pisa Cl. Sci. (5) \textbf{11}
  (2012), no.~2, 303--340.

\bibitem{met-spi3}
\bysame, \emph{Kernel estimates for some elliptic operators with unbounded
  coefficients}, DCDS-A \textbf{32} (2012), 2285--2299.

\bibitem{met-spi-tac}
G.~Metafune, C.~Spina, and C.~Tacelli, \emph{Elliptic operators with unbounded
  diffusion and drift coefficients in ${L}^p$ spaces}, Adv. Diff. Equat
  \textbf{19} (2012), no.~5-6, 473--526.

\bibitem{met-spi-tac2}
\bysame, \emph{On a class of elliptic operators with unbounded diffusion
  coefficients}, Evol. Equ. Control Theory \textbf{3} (2014), no.~4, 671--680.

\bibitem{wang}
F.~Y. Wang, \emph{Functional inequalities and spectrum estimates:the infinite
  measure case}, J. Funct. Anal. \textbf{194} (2002), 288--310.

\end{thebibliography}
\bibliographystyle{amsplain}

\end{document}